\documentclass[11pt]{article}  
\usepackage{amsmath,amssymb}
\newtheorem{Theorem}{Theorem}[section]
\newtheorem{Lemma}[Theorem]{Lemma}
\newtheorem{example}[Theorem]{Example}

\newtheorem{remark}[Theorem]{Remark}

\newenvironment{Remark}{\begin{remark}\rm}{\end{remark}}
\newtheorem{Definition}[Theorem]{Definition}
\newtheorem{Algorithm}[Theorem]{Algorithm}
\def\C{{{{\rm {\mbox{\small l}}} \kern -.50em {\rm
        C}}}}
\def\v{\vspace{1ex}}
  
\def\n{\noindent}
\def\dis{\displaystyle }
\def\pa{\partial}
\def\ve{\varepsilon}
\def\qed{\mbox{\,\rule{.5em}{.5em}}}
\newenvironment{Proof}
{\bigskip\noindent {\bf Proof.}~\rm}{\bigskip}
\begin{document}
\title{\huge All Zeros of the Riemann Zeta Function in the Critical
  Strip are Located on the Critical Line and are Simple}

\author{{Frank Stenger}
\\{\sl University of Utah}
\\{\sl Salt Lake City, UT \ \ 84112}
\\{\sl email: sinc\_f@msn.com}}

\date{}

\thispagestyle{empty}

\maketitle

\v \n \begin{abstract}
  \v \n In this paper we study the function $G(z) := 
  \int_0^\infty y^{z-1}(1 + \exp(y))^{-1} \, dy$\,, for $z \in {\mathbf 
  C}$\,.  We derive a functional equation that relates
    $G(z)$ and  $G(1 - z)$ for all $z \in {\mathbf C}$\,, and we prove:
  \begin{itemize}
  \item That $G$ and the Riemann Zeta function $\zeta$ have
    exactly the same zeros in the {\it critical region} $D := \{z \in
    {\mathbf C}\,: \Re\,z \in (0,1)\}$\,;   
  \item The Riemann hypothesis, i.e., that all of the zeros
    of $G$ in $D$ are located on the {\em critical line} $:= \{z \in D
    : \Re\,z = 1/2\}$\,; and that
  \item All the zeros of the Riemann Zeta function located on the
    critical line are simple.
  \end{itemize} 
\end{abstract}

\v \n {\bf Keywords}: Riemann hypothesis, Fourier transforms, Schwarz
  reflection principle, Cauchy--Riemann equations 

\section{Introduction and Summary}
\setcounter{equation}{0}

\v \n The proof of the Riemann hypothesis is a problem that many
mathematicians consider to be the most important problem of
mathematics. Indeed, it is one of at most {\it seven} mathematics
problems for which the Clay Institute has offered a million dollars
for its solution.  To this end, the .pdf publication [1] of
Bombieri (see also {[2]}) presents an excellent summary -- along with
references, to papers and books, to connections with prime numbers, to
Fermat's last theorem, and to the work of authors who have shown that
the first 1.5 billion zeros of the Zeta function listed with
increasing imaginary parts are all simple -- all of which are related
to the mathematics of this subject.  Similarly {\it Wikipedia} of the
web [3] offers an excellent summary along with references about
this subject.  The magazine {\em Nature} recently published a related
article about a discovery by Y. Zhang, of a conjecture (see [4]) on
the spacing of prime numbers $p_n$ with increasing size. Physicists
have also published results on the Riemann hypothesis: in [5],
Meulens compares data about the Riemann hypothesis with solutions of
two dimensional Navier Stokes equations, while others [6] have
compared eigenvalues of self--adjoint operators with zeros of the
Riemann Zeta function.  Several papers about solutions to the Riemann
hypothesis have also appeared.  To this end, the the papers of Violi
[7], Coranson--Beadu, [8], Garcia--Morales, [9], and Chen
[10] are similar to ours, in that their proof of the Riemann
hypothesis are for functions that are different from the Zeta
function, but which have the same zeros in $D$ as the Zeta
function. We note the incredible amount of numerical work that has
occurred, to find a zero of the zeta function in the critical strip
but not on the critical line: see e.g. [11,12,13].  We also
wish to mention the excellent work of [14], to establish that that all
zeros of the Riemann zeta function are on the critical line, provided
that all the zeros of the zeta function are simple.  To this end, our
present paper appears to be the only one which proves that all of the
zeros of the Zeta function on the critical line are simple.

\v \n Castelvecchi, author of the article [4] makes the comment:
``The Riemann hypothesis will probably remain at the top of
mathematicians wish lists for many years to come.  Despite its
importance, no attempts so far have made much progress.''

\v \n We wish of course to disagree with Castelvecchi's comment at the
end of the previous paragraph, since we believe that we have indeed
proved the Riemann hypothesis in this self--contained paper, in which
we accomplish the following:

\begin{enumerate}
\item In \S 2, defining the function $G$ and showing in
  detail that $G$ has exactly the same zeros, in the {\it critical
   strip}, $D := \{z \in {\mathbf C}: 0 < \Re\,z < 1\}$\,, including
  multiplicity, as the Zeta function;
\item Deriving a functional relation for $G$ in \S 2;  
\item Proving the positivity of $\Re\,G^{[2\,m]}(\sigma)$ for $\sigma
  \in (0,1)$ and the negativity of $\Im\,G^{[2\,m+1]}(\sigma)$ for
  $\sigma \in (0,1/2]$ in \S 2;  where 
  
    \begin{equation} 
    \begin{array}{c}
    \Re\,G^{[2\,m]}(\sigma + i\,t) := 
      \Re\,\left(\frac{\pa }{\pa 
      \sigma}\right)^{2\,m}\,\Re\,G(\sigma + i\,t) \\
    \\ 
   \Re\,(-\,i\,G^{[2\,m+1]}(\sigma + i\,t)) := 
     \Re\,(-\,i\,\left(\frac{\pa }{\pa
      \sigma}\right)^{2\,m+1}\,G(\sigma + i\,t)\,; 
    \end{array}
    \label{equation:1.1} 
    \end{equation} 
    
\item Introducing the Schwarz reflection principle in
  \S 3, which the functions $G^{[2\,m]}$ and $-\,i\,G^{[2\,m+1]}$ of
 Equation (\ref{equation:1.1}) satisfy in $D$\,; 
\item Proving by contradiction, via use of results developed in \S2
  and in \S4, that all of the zeros of G, (i.e., all of the zeros of
  the Zeta function) on the critical line $L$ are simple; and 
\item In \S 4, proving the Riemann hypothesis by contradiction, by use
  of results developed in \S 2, and in \S 3, of this paper, i.e., 
  proving that $G$ (and the zeta function) have no zeros in the region
  $D \setminus L$\,, where $L$ denotes the critical line, $L := \{z
  \in {\mathbf C}: \Re\,z = 1/2\}$\,.
\end{enumerate}

\v \n Let ${\mathbf R}^+ := (0,\infty)$\,, and let ${\mathbf R}^- =
(-\,\infty,0)$\,.  In this paper we thus derive results about the
function $G$ defined by the integral,

\begin{equation} 
  G(z) := \dis \int_{{\mathbf R}^+} \frac{y^{z-1}}{e^y + 1}\,dy\,,
  \ \ \Re\,z > 0\,,
  \label{equation:1.2}
\end{equation}

\n which is related to the well--known integral for the Riemann Zeta
function, where 

\begin{equation} 
  \zeta(z) := \dis \frac{1}{\Gamma(z)}\,\int_{{\mathbf R}^+}
  \frac{y^{z-1}}{e^y - 1}\,dy\,, \ \ \Re\,z > 1\,, 
  \label{equation:1.3}
\end{equation}

\n and where $\Gamma$ denotes the gamma function.

\v \n The operations of Schwarz reflection, the evaluation of $\Re\,G$
and $\Im\,G$ on important intervals of ${\mathbf R}$\,, can be readily
obtained from the Fourier transform representation of $G$\,, which is
gotten from Equation (\ref{equation:1.2}) and which is defined for
$\Re\,z > 0$\,, whereas an explicit Fourier transform of $\zeta$
defined by Equation (\ref{equation:1.3}) for $\Re\,z > 1$ does not
seem to be available.

\section{Representations and related properties of $G$}
\setcounter{equation}{0}

 \v \n In this section, we introduce the Fourier transform 
representation of $G$\,, and we develop related properties of $G$
which we shall require in our proof in \S 4.

\v \n The function $\zeta$ has many related properties,
with the best known of these given by:

\begin{equation} 
  \begin{array}{rcl}
    \zeta(z) & := & \dis \sum_{n=1}^\infty \frac{1}{n^z}\,\ \ \Re z >
       1\,, \\
  && \\
  \eta(z) & := & \dis \sum_{n=1}^\infty \frac{(-1)^{n-1}}{n^z}\,,\ \
       \Re z > 0\,,\ {\rm and} \\
  && \\
  G(z) & := &  \Gamma(z)\,\eta(z) := \left(1 -
  2^{1-z}\right)\,\Gamma(z)\,\zeta(z)\,,\ z \in {\mathbf C}\,. 
  \end{array}
  \label{equation:2.1}
\end{equation}

\n By setting $y = e^x$ and $z = \sigma + i\,t$ in
(\ref{equation:2.1}), we get a {\it Fourier integral representation}
of $G$\,, namely,

\begin{equation} 
 G(\sigma + i\,t) :=  \dis \int_{\mathbf R} \kappa(\sigma,x) \,
 e^{i\,x\,t} \, dx\,,
 \label{equation:2.2}
\end{equation}

\n where $\sigma \in (0,1)$\,, $t \in {\mathbf R}$\,, and where $\kappa$
is defined by  

\begin{equation}
  \kappa(\sigma,x) :=  \dis \frac{e^{\sigma\,x}}{1 + \exp(e^x)}\,.
  \label{equation:2.3}
\end{equation}

\subsection{Properties of $\kappa$\,, $\zeta$ and $G$} 

\v \n In this section we use the definition of $G$ given in Equation
(\ref{equation:1.2}) and the identities of Equation
(\ref{equation:2.1}) to derive a functional equation for $G$\,, and to
derive additional properties of $\kappa$ and $G$\,.  We also show in
detail, that $\zeta$ and $G$ have exactly the same zeros in $D$\,,
including multiplicity, that $\kappa(\sigma,-\,x) - \kappa(\sigma,x)$
is positive for all $(\sigma,x) \in (0,1/2] \times {\mathbf R}^+$\,, 
and strictly decreasing as a function of $\sigma$\,, for $\sigma \in
(0,1/2]$\,, and we determine ranges of values of $\Re\,G^{(m)}$\,,
$\Im\,G^{(m)}$ and their derivatives on the real line.
 
\v \n Let us next assign notations for the left and right half
of the complex plane, the critical strip(s), and the critical line.

\begin{Definition}
  \label{definition:2.1}
  \v \n Let ${\mathbf C}^-$ denote the left half of the complex plane,
  i.e., ${\mathbf C}^- := \{z \in {\mathbf C}\,:\,\Re\,z < 0\}$\,, and
  let ${\mathbf C}^+ := \{z \in {\mathbf C}\,: \Re\,z > 0\}$ denote
  the right half of the complex plane.  Let the critical strip be
  defined by $D = \{z \in {\mathbf C}: 0 < \Re\,z < 1\}$\,, and let
  the negative and positive critical strips be defined by $D^- = \{z
  \in D: \Im\,z \leq 0\}$\,, and $D^+ := \{z \in D\,: \Im\,z \geq 0\}$\,.
  The critical line is defined by $L := \{z \in D\,: \Re\,z =
  1/2\}$\,, and the negative and positive critical lines are defined
  by $L^- = \{z \in L : \Im\,z < 0\}$ and $L^+ = \{z \in L : \Im\,z >
  0\}$.   Let $a$ and $b$ denote two distinct points of $D$\,, let
  $\ell(a,b)$ denote the open line segment joining $a$ and $b$\,, and
  similarly for the half open, or closed, line segments joining $a$ and
  $b$\,.
\end{Definition}

\subsection{Relevant Gamma function relations} 

\v \n We shall require the use of the following lemma:

\begin{Lemma} 
  \label{lemma:2.2}
{(i.)} By replacing $z$ with $z/2$ in the duplication formula for the
Gamma function, we get:

$$
\Gamma(z) = (2\,\pi)^{-\,1/2}\,2^{z-1/2}\,\Gamma((z +
  1)/2)\,\Gamma(z/2)\,; \\
$$

\n {(ii.)} Both $\Gamma(1/2 + i\,x)$ and $\Gamma(1 + i\,x)$
   are bounded by $\pi^{1/2}$ for all $x \in {\mathbf R}$\,, by
   Equations (6.1.30) and (6.1.31) of [15]; and \\
\n {(iii.)} The function $1/\Gamma(z)$ is an entire function; $\Gamma$ 
   is analytic in ${\mathbf C}$ except for simple poles  
   at $z = -\,n$ ($n = 0\,,\ 1\,,\ 2\,,\ \ldots $).
\end{Lemma}

\begin{Proof} 
  Item (i.) Is just Equation (6.1.18) of [15] with z replaced by
  $z/2$\,; \\
  Items (ii.) Follow from Equations (6.1.30) and (6.1.31)
  of [15]; and \\
  Item (iii.) Is just a restatement of a result found
  in Chapter 16. of [15]. \end{Proof} \qed 

\subsection{Bounds on  $\kappa$} 
\v \n The next lemma describes some asymptotic bounds on the function
$\kappa$\,, which are obtained by inspection of Equation
(\ref{equation:2.3}).

\begin{Lemma}
  \label{lemma:2.3}
 For any $\ve \in (0,\sigma) \subset (0,1)$ and for $x$ real, we have

 \begin{equation}
   \kappa(\sigma,x) =  \left\{\begin{array}{l} {\mathcal
     O}\left(e^{(\sigma - \ve)\,x}\right)\,,\ \ x \to - \infty\,, \\
   {\mathcal O}\left(\exp\left(\sigma\,x - e^{x\,(1 -
     \ve)}\right)\,\right)\ \ x \to \infty\,.
    \end{array}
   \right.
   \label{equation:2.4} 
 \end{equation}

 \v \n Hence the integral $\int_{\mathbf R}
 Q(x)\,\kappa(\sigma,x)\,dx$ is finite for any polynomial $Q$\,.
\end{Lemma}

\begin{Proof}  The bounds of $\kappa$ given in
  (\ref{equation:2.4}) follow by inspection of the function $\kappa$
  as defined in (\ref{equation:2.3})\,. 
\end{Proof}\qed

\subsection{Analyticity definition of multiplicity} 

\v \n \begin{Definition}
  \label{definition:2.4}
  Let $z_0 \in {\mathbf C}$\,, let $m$ denote an integer, and let $f$
  be analytic in a neighborhood of $z_0$\,. 

  (a.) The function $f$ is said to have multiplicity $m$ at $z_0$ if
  there exists a finite number $c$ such that 

  $$
  \dis \lim_{z \to z_0} \frac{f(z)}{(z - z_0)^m} = c\,;
  $$

  (b.) If the multiplicity of $f$ at $z_0$ is $m$\,, and if $c \neq
  0$\,, then we shall more specifically say that $f$ is of exact
  multiplicity $m$ at $z_0$\,;

  (c.) If $f$ is of exact multiplicity $m$ at $z_0$\,, then $z_0$ is
  said to be a zero (resp., a pole) of $f$ of multiplicity $m$ if $m >
  0$ (resp., if $m < 0$).  In particular, if $m = 1$\,, (resp., if $m
  = - 1$\,,) then $z_0$ is said to be a simple zero (resp., a simple
  pole) of $f$\,;

  (d.) We shall say that $f$ vanishes with total multiplicity $m$ on a
  line segment $\ell(a,b)$ (or on $\ell(a,b]$\,, $\ell[a,b)$ or
  $\ell[a,b]$), if there exist a total of $n$ distinct points
  $z_1\,,\ z_2\,,\ \ldots\,,\ z_n \in \ell(a,b)$ (or on $\ell(a,b]$\,,
  $\ell[a,b)$ or $\ell[a,b]$), with $z_j$ of multiplicity
  $\mu_j \geq 0$\,, for $j = 1\,,\ 2\,,\ \ldots, n$\,, such that

  $$
  m := \dis \sum_{j=1}^n \mu_j \,, 
  $$

  \n and we shall say that this total multiplicity on $\ell[a,b]$ is
  exact if $f^{(\mu_j)}(z_j) \neq 0$ for all $z_j =
  z_1\,,\ z_2\,,\ \ldots\,,\ z_n \in \ell[a,b]$\,. 
  
 \end{Definition}

\subsection{Functional equations for $\zeta$ and $G$} 

\v \n An important identity of the Riemann Zeta function is the well
known functional equation for $\zeta$, due to Riemann: 

\begin{equation}
   \pi^{-\,(1-z)/2}\,\Gamma((1 - z)/2)\,\zeta(1 - z) = \pi^{-\,z/2}\,
  \Gamma(z/2)\,\zeta(z)\,.
    \label{equation:2.5}
 \end{equation}  

 \v \n This functional equation for the Riemann Zeta function has many
 important uses, including, e.g., the analytic continuation of Zeta
 to all of ${\mathbf C}$\,.

 \v \n The function $G$ also possesses a functional equation which is
 given in Lemma \ref{lemma:2.5} below, which plays a similar role as
 the functional equation for $\zeta$\,.  The functional equation for
 $G$ is gotten by substituting the right--hand--side of the third
 equation of (\ref{equation:2.1}) into (\ref{equation:2.5}), and by
 use of Lemma \ref{lemma:2.2}:

\begin{Lemma}
\label{lemma:2.5}
Let $z \in {\mathbf C}$\,, and let $G$ be defined as in
(\ref{equation:2.2}).  Then, a functional equation for the function
$G$, valid for all $z \in {\mathbf C}$ is:

\begin{equation} 
  \label{equation:2.6}
  \frac{2^{1\,-\,z} -
  1}{(4\,\pi)^{\frac{1-z}{2}}\,\Gamma\left(\frac{1 + (1 
      - z)}{2}\right)}\,G(1 - z) :=
  \frac{2^z - 1}{(4\,\pi)^{\frac{z}{2}}\,\Gamma\left(\frac{1 +
    z}{2}\right)}\,G(z)\,. 
\end{equation}

\n This equation can also be written in a more compact form: 

\begin{equation} 
  G(1 - z) := K(z)\,G(z)\,,
    \label{equation:2.7}
\end{equation}

\n where 

\begin{equation} 
  \begin{array}{rcl}
  K(z) & := & \dis (4\,\pi)^{1/2-z}\ \frac{2^z - 1}{2^{1-z} - 1}
  \frac{\Gamma\left(\frac{1 + (1 - z)}{2}\right)}{\Gamma\left(\frac{1 +
      z}{2}\right)}\,, \\
  && \\
  G^{[m]}(1 - z) & = & \dis (-1)^m\,\sum_{j=0}^m {m \choose j}
  K^{[m-j]}(z)\,G^{[j]}(z)\,.
  \end{array}
   \label{equation:2.8}
\end{equation}

\n and where $K$ is non--vanishing in $D$\,. 
\end{Lemma}

\begin{Proof}
The identity (\ref{equation:2.6}) is valid for all $z \in {\mathbf C}$
since equation (\ref{equation:2.5}) is valid for all $z \in {\mathbf
  C}$\,, and since, except for simple poles, all of the functions
multiplying the zeta function in the third equation of
(\ref{equation:2.1}) are analytic in all of ${\mathbf C}$\,.  An
alternate proof obtains by inspection of Equations
(\ref{equation:2.7}) and (\ref{equation:2.8}).  That $K$ is analytic
and non--vanishing in $D$ follows from Lemma \ref{lemma:2.2} and by
inspection of Equation (\ref{equation:2.8}).  \end{Proof} \qed

\subsection{Zeros of $G$ and $\zeta$ in $D$} 

\v \n We prove here the $G$ and $\zeta$ have the same zeros with the
same multiplicity in $D$\,, and that these zeros are isolated.  

 \begin{Lemma}
  \label{lemma:2.6}
  (i.) The functions $G$ and $\zeta$ have exactly the same zeros in
  $D$\,, including multiplicity; and \\
  (ii.) All zeros of $G$ in ${\mathbf C}$ are isolated. 
 \end{Lemma}
 
\begin{Proof} 
  (i.) By inspection if the third equation of (\ref{equation:2.1}) we
  get, if $z_0 \in D$ is a zero of $G$ of multiplicity $k \geq 1$\,,
  then $(*)\,G(z) = w(z)\,\zeta(z)$\,, where the
  function $w(z) := (1 - 2^{1-z})\,\Gamma(z)$ is analytic and
  non--vanishing in $D$\,, so that $z_0$ is also a zero of $\zeta$ of
  multiplicity $k \geq 1$\,.  In addition, by taking the $n^{th}$
  derivative of both sides of $(*)$\,, with non--negative integer 
  $n$\,, we get,

  \begin{equation} 
    G^{(n)}(z) = \dis (-1)^n\,\sum_{j=0}^n {n \choose j}
    w^{(n-j)}(z)\,\zeta^{(j)}(z)\,.   
  \label{equation:2.9}
  \end{equation}

  \n Hence, if $z_0 \in \Omega^+$ is a zero of $\zeta$ of multiplicity
  $m \geq 1$\,, then by applying induction to Equation
  (\ref{equation:2.9}), with respect to $n =
  0\,,\ 1\,,\ 2\,,\ \ldots\,, m$ we conclude that the multiplicity of
  the zero $z_0$ of $G$ is also $m$\,;

  \v \n (ii.)  Suppose that there exists a cluster of zeros
  $\{z_j\}_{j=1}^\infty$ of $G$ in $D$ with a sub-sequence that has a
  limit point $z$.  If $z \in D$\,, then $G$ would have to vanish, by
  Vitali's theorem.  If $z$ is on the line $\{\Re\,z = 0\}$\,,  then,
  since $D \subset {\mathbf C}^+$\,, and since $G$ is analytic in
  ${\mathbf C}^+$\,, it follows by use of the functional equation of
  $G$, that $0 = G(z) = G(1 - z)$ where the point $1 - z$ is now
  located on the line $\{z \in {\mathbf C} : \Re\,z = 1\,\}$\,, i.e.,
  we are back to the previous case of the convergence of such a
  sub-sequence to a point on the interior of the right half plane,
  where $G(z)$ is analytic and bounded, so that $G(z)$ would again
  have to vanish identically in ${\mathbf C}$\,.
  \end{Proof} \qed

\subsection{Definitions of $\kappa^\mp\,,\ G^{(m)}$ and $G^{[m]}$}

\begin{Definition}
  \label{definition:2.7}
  Let $G$ and $\kappa$ be defined as in Equation (\ref{equation:2.2}),
  and let us define $\kappa^\mp(\sigma,x)$ as follows:

   \begin{equation}
       \kappa^\mp(\sigma,x) := \left\{
      \begin{array}{l}
        \kappa(\sigma,\mp x)\,,\ x  \in {\mathbf R}^+\,,\\
        0\,,\  x  \in  {\mathbf R}^-\,,
      \end{array}
      \right.\\
      \label{equation:2.10}
   \end{equation}

   \v \n If for brevity, we write $\kappa^\mp$ for
   $\kappa^\mp(\sigma,x)\,,\ C$ and $S$ for $\cos(x\,t)$ and
   $\sin(x\,t)$\,, and $\int$ for $\int_{{\mathbf R}^+}$\,, then 
   Equation (\ref{equation:1.2}) yields the following definitions for
   $G^{(2\,n)}$ and for $G^{(2\,n+1)}$\,, where $n$ denotes a
   non--negative integer: 

   \begin{equation} 
      \begin{array}{l}
      G^{(2\,n)}(\sigma\,+\,i\,t)   \\
    \\ \ \ =  \dis
    (-1)^n\,\int\,x^{2\,n}\,\left((\kappa^- + \kappa^+)\,C -
    i\,(\kappa^- - \kappa^+)\,S\right)\,dx\,, \\
    \\
    G^{(2\,n+1)}(\sigma\,+\,i\,t) \\
    \\ \ \ = \dis 
    (-1)^{n+1}\,\int\,x^{2\,n+1}\,\left((\kappa^- + \kappa^+)\,S
    + i\,(\kappa^- - \kappa^+)\,C\,\right)\,dx\,. 
    \end{array}
      \label{equation:2.11}
    \end{equation}

    \n Let us define $G^{[m]}$ for any non--negative integer $m$ by

    \begin{equation} 
    G^{[m]}(\sigma + i\,t) =
    \left(\frac{\pa}{\pa\,\sigma}\right)^m\,G(\sigma + i\,t)\,,
    \label{equation:2.12}
    \end{equation}

    \n Then, 

    \begin{equation} 
      G^{(m)}(\sigma + i\,t)\, = \,  
      \dis i^m\,\left(\frac{\pa}{\pa\,(i\,t)}\right)^m\,G(\sigma +
      i\,t) \, = \, 
      \dis i^m\,G^{[m]}(\sigma + i\,t)\,,
      \label{equation:2.13}
    \end{equation}
    
\n where these functions are readily shown to exist, by Lemma
\ref{lemma:2.3}. 
\end{Definition}

\begin{Lemma}
  \label{lemma:2.8}
Let the functions $G^{[m]}$ be defined as in Definition
\ref{definition:2.7}\,.  Then, for all $m =
0\,,\ 1\,,\ 2\,,\ \ldots$\,, and for all $\sigma \in (0,1)$\,,
$G^{[m]}$ is analytic on the right half plane, and hence also in
$D$\,.  In particular given any $\ve > 0$\,, $G^{[m]}(z)$ is uniformly
bounded in the region $\{z \in D\,: \Re\,z \geq \ve\}$\,.
\end{Lemma}

\begin{Proof} This result follows directly by inspection of Equation
  (\ref{equation:2.2}) and Lemma \ref{lemma:2.3}.  We omit the
  straight--forward proofs. \end{Proof} \qed

\subsection{Restricting the domain of $\Im\,G^{(2\,m+1)}$} 

\v \n The following lemma restricts the domain of some
of our identities:

\begin{Lemma} 
  \label{lemma:2.9}
  Let $\Delta$ be defined by $\Delta(\sigma,x) := \kappa^-(\sigma,x) -
  \kappa^+(\sigma,x)$\,, where the functions $\kappa^\mp$ are defined
  in Definition \ref{definition:2.7}.  Then $\Delta(\sigma,x) > 0$
  for all $(\sigma,x) \in (0,1/2] \times {\mathbf R}^+$\,, and
  moreover, $\Delta(\sigma,x)$ is a strictly decreasing function of
  $\sigma \in (0,1/2]$\, for every fixed $x \in {\mathbf R}^+$\,.
\end{Lemma}

\begin{Proof} We have 

  $$
  \frac{\pa}{\pa \sigma}\,\Delta(\sigma,x) = -\,x\,(\kappa^-(\sigma,x) +
  \kappa^+(\sigma,x))\,,
  $$

  \n which shows $\Delta(\sigma,x)$ is a strictly decreasing function
  of $\sigma \in (0,1/2]$ for all fixed $x \in {\mathbf R}^+$\,.  By
  making the one--to--one transformation $y = e^{-\,x/2}$ of
  ${\mathbf R}^+ \to (0,1)$ in the above expressions for
  $\Delta(\sigma,x)$\,, and then setting $W := (1 + e^{y^2})(1 +
  e^{1/y^2})$\,, we get 

  \begin{equation} 
  \begin{array}{rcl}
    \Delta(\sigma,x) & := &  \delta(\sigma,y)/W\,, \\
    && \\
    \delta(\sigma,y) & := & y^{2\,\sigma}\,(1 + \exp(1/y^2)) -
    1/y^{2\,\sigma}\,(1 + \exp(y^2))\,.
  \end{array}
  \label{equation:2.14}
  \end{equation}

  \n Since $1/W$ is positive on $(0,1)$\,, and
  since $\Delta(\sigma,x)$ is a strictly decreasing function of
  $\sigma \in (0,1/2]$ for all fixed $x \in {\mathbf R}^+$\,, we need
  only prove that $\delta(1/2,y) > 0$ for all $y \in (0,1)$\,.  To 
  this end we have, by use of Taylor series expansions, that 

 \begin{equation}
   \begin{array}{rcl}
     \delta(1/2,y) & = & \dis y\,(1 + e^{1/y^2})\,-\,1/y\,(1 +
     e^{y^2})\,,\ \ y \in (0,1) \\   
     && \\
     & = & \dis \,(1/2)(1/y - y)\,(1/y^2 + 1 + y^2 - 2) \\
     \\
     && \ + \dis \sum_{n=3}^\infty \frac{1/y^{2\,n-1} -
       y^{2\,n-1}}{n\,!} > 0\,,\ \ y \in (0,1)\,.
   \end{array}
  \label{equation:2.15}  
\end{equation}
 
\n Upon proceeding from the first to the second line of Equation
(\ref{equation:2.15}) we used the following relations, which are valid
for all $y \in (0,1): 0 < y < 1 < 1/y$\,.  The right--hand side of
Equation (\ref{equation:2.15}) then shows that $\delta(1/2,y) > 0$ for
all $y \in (0,1)$\,, i.e., that $\Delta(\sigma,x) > 0$ for all
$(\sigma,x) \in (0,1/2] \times {\mathbf R}^+$\,. \end{Proof} \qed

\subsection{Inequalities for $G^{(m)}(\sigma)$} 

\v \n The following lemma summarizes values of $G^{(m)}(\sigma)$
that have been established above.

\newpage
\begin{Lemma}
  \label{lemma:2.10}
  Let $m$ denote a non--negative integer. Then:
\begin{description}  
\item{(i.)} $(-1)^m\,\Re\,G^{(2\,m)}(\sigma) = \Re\,G^{[2\,m]}(\sigma)
  > 0$ for all $\sigma \in (0,1)$\,;
\item{(ii.)} $(-1)^m\,\Im\,G^{(2\,m+1)}(\sigma) =
\Re\,G^{[2\,m+1]}(\sigma) < 0$ for all $\sigma \in (0,1/2]$\,; and  
\item{(iii.)} $\Im\,G^{(2\,m)}(\sigma) = \Im\,G^{[2\,m]}(\sigma) 
  = \Im\,G^{[2\,m+1]}(\sigma) = \Re\,G^{(2\,m+1)}(\sigma) = 0$   
  for all $\sigma \in (0,1)$\,.
\end{description}
\end{Lemma} 

\begin{Proof} 
  The proof of each Item of Lemma \ref{lemma:2.10} follows by
  inspection of Equation (\ref{equation:2.11}).  \end{Proof} \qed

\begin{Remark} 
  \label{remark:2.11}
  \begin{enumerate}
  \item The functions $\Re\,G^{(2\,m)}(z)$\,, $\Re\,G^{[2\,m]}(z)$\,,
    $\Im\,G^{(2\,m+1)}(z)$\,, and  $\Re\,G^{[2\,m+1]}(z)$ 
    yield the correct sign of of these functions at $z = 1/2$\,;
  \item Not so for the functions $\Im\,G^{(2\,m)}(z)$\,,
    $\Im\,G^{[2\,m]}(z)$\,, $\Im\,G^{[2\,m+1]}(z)$ or
    $\Re\,G^{(2\,m+1)}(z)$), which vanish identically on
    $\ell((0,1)$\,; and 
  \item If any of the functions of the previous item vanish on
    $\ell(0,1)$\,, then it will be better to use $G^{(2\,m)}$\,,
    $G^{[2\,m]}$\,, $G^{(2\,m+1)}$\,, or $G^{[2\,m+1])}$\,, in place 
    of them.
\end{enumerate}
\end{Remark}

\section{Schwarz reflection}
\setcounter{equation}{0}

\v \n We present the Schwarz reflection principle, which we define
as follows:

\begin{Definition}
  \label{definition:3.1}
  \v \n Let $f$ be analytic in $D$\,, and real on $(0,a)$\,, for
  some $a \in [1/2,1)$\,.  Then $f$ can be continued analytically
  (i.e., reflected) across $(0,a)$ from $D^\mp$ to $D^\pm$) by
  means of the formula
  
  \begin{equation} 
    f(\overline{z}) = \overline{f(z)}\,.
    \label{equation:3.1}
  \end{equation}

\end{Definition}

  \begin{Remark}
  \label{remark:3.2}
  The Schwarz reflection principle enables analytic continuation from 
  $D^\mp$ to all of $D$\,.  For example, if $n$ denotes a
  non--negative integer, so that the functions $f(z^+) :=
  G^{(2\,n)}(z^+)$ and $g(z^+) := -\,i\,G^{(2\,n+1)}(z^+)$ are given
  for $z^+ = \sigma + i\,t^+ \in D^+$\,, then by Lemma
  \ref{lemma:2.10}, and by Equation (\ref{equation:2.8}),
  $\Re\,f(\sigma)$ is a non--vanishing function of $\sigma$ on
  $(0,1)$\,, while if $\sigma \in (0,1/2]$\,, then $\Re\,G(\sigma)$ is 
  a non--vanishing function of $\sigma$\,, $\Im\,G(\sigma + i\,t)$
  changes sign as $t$ changes sign, while $\Re\,G(\sigma + i\,t$ does 
  not change sign as $t$ changes sign. 
\end{Remark}

\section{Proof of the Riemann Hypothesis}
\setcounter{equation}{0}

\v \n This section contains the proof of the following theorem: 

\begin{Theorem}
  \label{theorem:4.1}
  Let $D := \{z \in {\mathbf C}: \Re\,z \in (0,1)\}$ denote the
  critical region in the complex plane, and let $L = \{z \in D: \Re\,z
  = 1/2\}$\,.  Every zero $z$ of $G$ located on the critical line $L$ 
  is a simple zero of $G$\,, and every zero $z$ of $G$ in $D$ is 
  is located on the critical line $L$\,. 
\end{Theorem} {

\begin{Proof} We shall prove this theorem by means of the proof by
induction, and by contradiction, in {\it two} lemmas, which follow. 
We are able to carry out our proofs using $G$ and it's derivatives
instead of those of $\Re\,G$ or $\Im\,G$\,, see Remark
\ref{remark:2.11}.

\v \n The zeros of $G$ on $L$ are isolated, by Lemma \ref{lemma:2.6}.
Hardy and Littlewood proved in [16] that the Zeta function, and
hence $G$\,, has an infinite number of zeros on the critical line
$L$\,.  We begin this section by proving that all of the zeros of $G$
on $L$ are simple.

\begin{Lemma} 
  \label{lemma:4.2}
Every zero $\xi = 1/2 + i\,\tau$ of $G$ on the critical line $L :=
\{\xi \in {\mathbf C}: \Re\,\xi = 1/2\}$ is a simple zero of $G$\,.
\end{Lemma}

\begin{Proof}  Let us label all of the zeros $\xi_j = 1/2 +
i\,\tau_j$ of $G$ on $L^+$ in increasing order of imaginary parts as
follows: $\xi_j = 1/2 + i\,\tau_j$\,, where $0 < \tau_1 < \tau_2 <
\ldots$\,.

\v \n We wish to prove that {\it the case of a zero $\xi_j \in L$ of
$G$ of exact multiplicity greater than {\it one} is not allowed.}
Hence: {\it Let $P(m)$ be the proposition that the zeros of $G$ on $L$
are simple, for $m = 1\,,\ 2\,,\ \ldots\,,$.}  

\v \n {\it (I.) Proof that the case of a zero $\xi_1$ of $G$ of exact
multiplicity greater than {\it one} is not allowed.}  To this end,
let us assume the contrary, i.e., that each of $G(\xi_1)$ and
$G(\overline{\xi_1}$ vanish with multiplicity $\nu_1 + 2$\,, where
$\nu_1$ denotes a non--negative integer.  Then $G$ vanishes on
$\ell[\xi_1,\overline{\xi_1}]$ with total multiplicity $2\,\nu_1 +
4$\,, or equivalently, $G^{(2\,\nu_1 + 4)}$ vanishes on
$\ell[\xi_1,\overline{\xi_1}]$ with multiplicity {\it one}, which is
not allowed, by Lemma \ref{lemma:2.10}.  This proves that the zero
$\xi_1$ of $G$ of exact multiplicity greater than {\it one} is not
allowed.

{\it (I.2) Proof that if $P(n)$ is true for all integers $n$ with 
$1 \leq n \leq k$\,, for some $k \geq 1$\,, then $P(n)$ is true for 
$n = k + 1$\,.}  If $P(n)$ is true for $1 \leq n \leq k$\,, then
each of $G(1/2 \pm i\,\tau_n$ vanish with multiplicity {\it one} for
all $n \in (1,k)$\,.

\v \n Let us then assume for $n = k + 1$\,, that both
$G(1/2 \pm i\,\tau_{k+1})$ vanish with exact multiplicity $\nu_{k+1} +
2$\,, where $\nu_{k+1}$ denotes a non--negative integer.  Since $G$
vanishes at $2\,k$ points on $\ell(\overline{\xi_{k+1}},\xi_{k+1})$\,,
$G$ vanishes on $\ell[\overline{\xi_{k+1}},\xi_{k+1}]$ with total
multiplicity $2\,\nu_{k+1} + 2\,k + 4$\,.  Equivalently,
$G^{(2\,\nu_{k+1}+2\,k+4)}(1/2)$ vanishes on
$\ell[\xi_{k+1},\overline{\xi_{k+1}}$ with multiplicity {\it one}, but
which is not allowed, by lemma \ref{lemma:2.10}.

\v \n This proves that the zero $\xi_{k+1} = 1/2 + i\,\tau_{k+1}$ of
$G$ of exact multiplicity greater than {\it one} is not allowed.

\v \n This completes the proof of Lemma \ref{lemma:4.2}. \end{Proof}
\qed

\begin{Lemma} 
  \label{lemma:4.3}
  If $z = s + i\,\tau \in D$\,, with $(s,\tau) \in \in (0,1) \times
  {\mathbf R}$\,, and with $s \neq 1/2$\,, then $G(z) \neq 0$\,.
\end{Lemma}

\v \n \begin{Proof} If $\tau = 0$\,, then $G(s + i\,\tau) = G(s) 
  \neq 0$ for all $s \in (0,1)$\,, by Lemma \ref{lemma:2.10}, so we
  shall assume, in the remainder of our proof which follows, that
  $\tau \neq 0$\,.

  \v \n We shall examine the vanishing of $G$ at points $\xi :=   
  \xi_{j,k} := s_{j,k} + i\,\tau_j \in D$\,, where for integer $m >
  0$\,, we have $0 < \tau_1 < \tau_2 < \ldots < \tau_m$\,, and where
  for each fixed $j \in (1,m)$\,, with $\mu_j$ a positive integer, we 
  shall assume that $0 < s_{j,1} < s_{j,2} < \ldots < s_{j,\mu_j} <
  1/2$\,.  For simplicity of proof, we shall consider, for
  each fixed $j \in (1,n)$, the appearance of $\xi_{j,k}$ in the order
  $k = \mu_j,\ \mu_j - 1,\ \ldots\,\ 1$\,.

  \v \n Let us state the induction hypothesis: {\it Let $m$ denote a
  positive integer, and let $P(m)$ be the proposition that
  $G(\xi_{j,k}) \neq 0$ for all $(j,k) \in 
  \{(1,m) \times (1,\mu_j)\}$\,.}
  
  \v \n {\it I.1: Proof that $P(1)$ is true, for $\xi_{1,k}\,,\ k \in
    (1,\mu_1)$\,.} 

  \v \n Let us assume the opposite of the above induction hypothesis,
  i.e., that each of the two values, $G(s_{1,\mu_1} \pm i\,\tau_1)$
  vanishes with multiplicity $\nu_1 + 1$\,, where $\nu_1$ denotes a
  non--negative integer.  Then, we get, by reflection, that each of
  $G(s_{1,\mu_1} \pm i\,\tau_1)$ vanishes with multiplicity $\nu_1 +
  1$\,, where $\nu_1$ denotes a non--negative integer. Hence if we
  assume in addition, that $G(1/2 + i\,\tau_1) \neq 0$\,, we get, by
  reflection, that $G$ vanishes on $\ell[s_{1,\mu_1} + 
  i\,\tau_1,s_{1,\mu_1} - i\,\tau_1]$ with total multiplicity $2\,\mu_1 +
  2$\,.  Equivalently, $G^{(2\,\nu_1+2)}(s_{1,\mu_1})$ vanishes with
  multiplicity {\it  one}, but which is not allowed by Lemma
  \ref{lemma:2.10}.
  
  \v \n On the other hand, if each of $G(1/2 \pm i\,\tau_1)$ vanishes
  with multiplicity {\it one}, and if each of $G(s_{1,\mu_1} \pm
  i\,\tau_1)$ vanishes with multiplicity $\nu_1 + 1$\,, where $\nu_1$
  denotes a non--negative integer, then $G$ vanishes on
  $\ell[\xi_{1,\mu_1},\overline{\xi_{1,\mu_1}}]$ with total
  multiplicity $2\,\nu_1 + 4$\,.  Equivalently,
  $G^{(2\,\nu_1+4)}(s_{1,\mu_1})$ then vanishes with multiplicity {\it
  one}, but which is also not allowed, by Lemma \ref{lemma:2.10}.
  
  \v \n Hence the vanishing of $G(s_{1,\mu_1} + i\,\tau_1$ with
  positive multiplicity $\nu_1 + 1$ is not allowed, regardless of
  whether or not $G(1/2 + i\,\tau_1)$ vanishes. By the same argument,
  the vanishing of $G(s_{1,k} + i\,\tau_1$ is not allowed for $k = \mu_1 -
  1\,,\ \mu_1 - 2\,, \ldots\,,\ 1$\,, regardless of whether or not
  $G(1/2 + i\,\tau_1)$ vanishes.  Hence, the vanishing of
  $G(s_{1,k} + i\,\tau_1$ is not allowed for all $(1,k) \in \{(1) 
  \times (1,\mu_j)$\,.

  \v \n Hence $P(1)$ is true.

  \v \n {\it I.2: Proof that if $P(n)$ is true, for $1 \leq n \leq
    m$\,, where $n$ denotes a positive integer, then $P(n+1)$ is
    true.}

  \v \n If $P(n)$ is true, there are no points $\xi_{j,k} \in D^+
  \setminus \{L^+\}$\,, with $k \in (1,\mu_j)$ for every fixed $j \in
  (1,n)$\,, such that $G(\xi_{j,k}) = G(s_{j,k} + i\,\tau_j) = 0$\,.
  Hence we need to prove, that if $P(n)$ is true, then $P(n+1)$ is 
  true, i.e., that $G(\xi_{n+1,\mu_k}) = G(s_{n+1,\mu_k} + i\,\tau_{n+1})
  \neq 0$\,, for all $k \in (1,n+1)$\,.  
   
  \v \n Let us assume the opposite, i.e., let us assume
  that both $G(s_{n+1,\mu_{n+1}} \pm i\,\tau_{n+1})$ vanish with
  multiplicity $\nu_{n+1} + 1$\,, where $\nu_{n+1}$ denotes a
  non--negative integer, and let us assume, in addition, that $G(1/2 +  
  i\,\tau_{n+1}) \neq 0$\,. Then $G$ vanishes on       
  $\ell[\xi_{n+1,\mu_{n+1}},\overline{\xi_{n+1,\mu_{n+1}}}]$ with total
  multiplicity $2\,\nu_{n+1} + 2$\,, where $\nu_{n+1}$ denotes a
  non--negative integer.  Equivalently,
  $G^{(2\,\nu_{n+1}+2)}(s_{n+1,\mu_{n+1}})$ vanishes with multiplicity
  {\it one}\,, but which is not allowed, by Lemma \ref{lemma:2.10}.

  \v \n On the other hand, if $G(1/2 + i\,\tau_{n+1})$ vanishes, then
  it vanishes with multiplicity {\it one}, so that $G$ vanishes on
  $\ell[\xi_{n+1,\mu_{n+1}},\overline{\xi_{n+1,\mu_{n+1}}}]$ with
  total multiplicity $2\,\nu_{n+1} + 4$\,, or equivalently, 
  $G^{(2\,\nu_{n+1}+4)}(s_{n+1,\mu_{n+1}}) = 0$\,, but which is not
  allowed, by Lemma \ref{lemma:2.10}.  
  Hence we cannot allow the vanishing of $G(\xi_{n+1,\mu_{n+1}})$ with
  multiplicity $\nu_{n+1} + 1 > 0$\,, where $s_{n+1,\mu_{n+1}}$ is the
  largest such value on $(0,1/2)$\,.

  \v \n Next, the second largest, $s_{n+1,\mu_{n+1}-1}$ becomes the
  largest, and we can eliminate this case with the same proof that we
  used above to eliminate the case of $\xi_{n+1,\mu_{n+1}}$\,, and so
  on, until we have eliminated all of the vanishings of
  $G(\xi_{n+1,k})$\,, for $k =
  \mu_{n+1},\ \mu_{n+1}-1,\ \ldots\,\ 2,\ 1$\,.  Hence $P(n + 1)$ is
  true.

  \v \n This completes the proof of Lemma \ref{lemma:4.3}. \end{Proof} 
  \qed  

  \v \n This completes the proof of Theorem \ref{theorem:4.1}.
  \end{Proof} \qed }
  
\v \n {\bf Acknowledgments:} 

\v \n The author is grateful to Drs. Paul Gauthier and
  Marc Stromberg for their corrections of the proofs related to Lemma
  \ref{lemma:2.5}. 

\section*{References} 
\begin{description} 
  \begin{sloppypar}
\item{[1]} E. Bombieri, {\it the Riemann
  Hypothesis -- official problem description}, Clay Mathematical
      Institute (2014). {1}
\item{[2]}
  www.claymath.org/sites/default/files/official\_problem\_description.pdf 
\item{[3]} https://en.wikipedia.org/wiki/riemann\_zeta\_function
\item{[4]} D. Castelvetcchi, {\it Mathematician claims Prime Number
    Riddle Breakthrough} Nature {\bf 611} 645--646 (2022).
\item{[5]} Meulens, R., {\em The Proof of the Riemann Hypothesis
    and an Application to Physics},  APM, Applied Mathematics {\bf
    10} (2019) 967-988.
\item{[6]} S. Endres, and F. Steiner, {\it the Berry--Keating
  operator on $L^2({\mathbf R})$ and on compact quantum graphs with
  general self--adjoint realizations}, J. of Physics, A. Mathematical
  and Theoretical, {\bf 43} 095204 {\bf DOI}
  10.1088/1751-8113/43/9/095204.
\item{[7]} R. Violi,  {\em All complex zeros of the Riemann Zeta
  Function are on the Critical Line},
  https://arxiv.org/pdf/2010.05335.pdf
\item{[8]} J.--M. Coranson--Beaudu, {\em A Speedy
    New Proof of the Riemann Hypothesis}, {\bf 10} (2021) 62-67.Comp.,
       {\bf 46} (1986) 667--681.
\item{[9]} Garcia--Morales, {\it On the non--trivial zeros of the
   Dirichlet eta function}, arXiv 2007.04317v1 (220). 
    (2022) 374--391.
\item{[10]} C. Chen, {\it Proof of Riemann Conjecture}, APM {\bf 12}
    (2022) 374--391.
\item{[11]} R.P. Brent, {\it On the zeros of the Riemann zeta function
     in the critical strip}, Math. Comp. {\bf 33} (148) (1979)
     1361--1372.
\item{[12]} A.M. Odlyzko, {\it The $10^{22} $ Zero of the Riemann Zeta
  Function}, Dynamical, Spectral, and Arithmetic Zeta Functions, San
  Antonio, Texas (1999) Contemporary Mathematics, V. 290, AMS,
  Providence RI, (2001), pp. 139--144.
\item{[13]} H.M. Edwards, {\it Riemann's Zeta Function}, Academic Press
  (1974).
\item{[14]} D.A. Cardon and S de Gaston Roberts, {\it An equivalence
     for the Riemann Hypothesis in terms of orthogonal polynomials,}
     J. Approx. Theory 138 (2006), no. 1, 54-64. 
\item{[15]} Abramowitz, A. \& Stegun, I.A., {\em Handbook of
    mathematical functions}, National Bureau of Standards Applied
  Mathematical Series, {\bf 55}, 1964. {18}
\item{[16]} G.H. Hardy and J.E. Littlewood, {\it Contributions to the
  theory of the Riemann Zeta-Function and the theory of the
  distribution of primes}, Acta Mathematica, {\bf 41} (1917) 119-196.
 \end{sloppypar}
\end{description}

\end{document}